# Dynamical Characterization of Fractal Objects: Determination of the Fine Fractal Topology Using the Energy Cover

Luiz Bevilacqua, Marcelo Miranda Barros, Gil Márcio A. Silva

**Abstract.** The foundation of the theory presented here has already been proved to be effective for the case of curves belonging to the Koch family. The present paper extends the investigation to more complex curves, namely randomly generated curves and the Weierstrass-Mandelbrot curve. The analysis is focused on numerical experiments. The results obtained with the numerical analysis allow advancing some interesting proposition concerning the fine fractal structure of plane curves. The method uses the dynamical response of appropriate harmonic oscillators built up according to the geometry of the fractal set. For the plane motion three fundamental periods are determined. Each period plotted against the characteristic length of the corresponding curve in a log-log scale estimates the fractal dimension. It is shown that each degree of freedom corresponds to a distinct energy covering associate to a certain topological property. Using this technique it is possible to distinguish the fine fractal structure of plane curves. The dynamic technique is proved to work also for the identification problem, that is, to determine the fractal characteristics embedded in some sample of a given fractal set. It is proposed the classification of fractal structures into at least two categories, perfectly fractal and partially fractal. The affine similarity can be precisely detected using this method.

## INTRODUCTION

The main purpose of this paper is to explore the notion of a multiple fractal characterization of plane curves. Indeed the dynamical determination of the fractal dimension of plane objects leads to three independent figures. They may or may not coincide. The reason for the differences comes from three independent elastic energy distributions, that can be assumed as generalized energy covers, induced by three selected excitation. That for simple harmonic oscillators on a plane it is always possible to excite three independent frequencies is obvious from the three eigenvalues corresponding to the oscillator stiffness matrix. We will see that it is possible to specify one particular dynamical dimension that is equivalent to the Hausdorff dimension. The other two provide additional informations which are equally important to define the complete fractal structure of the object. Let us recall that the characterization of the fractal Hausdorff dimension is associated to a well-defined cover. Indeed consider a fractal object Ж, and a countable set of balls $B_i$ with radius $r$ defined as the largest distance of two points belonging to the ball. Let this set be a cover of Ж, that is, Ж belongs to ∩$B_i$. The outer measure of the cover is less or equal to the sum of the diameters of all balls of the cover. The s-dimensional Hausdorff outer measure H(s,r) is the infimum of all measures raised to the power "s" of all such covers, that is:

$$H(s,r) = \inf\left(\sum_i (diameterB_i)^s\right)$$

Moreover there exists one and only one number D, defined as the Hausdorff dimension of the fractal object such that:

$$\lim_{r \to 0}(H(s,r)) \to \infty \quad \text{for all} \quad s > D$$
$$\lim_{r \to 0}(H(s,r)) \to 0 \quad \text{for all} \quad s < D$$

Let us now recall the basic notions of the dynamical dimension. For well determined plane curves it has been proved that it is possible to recover the Hausdorff dimension using a new technique that we call the dynamical characterization of fractal objects[1]. Let us review the main characteristics of the method.

Given a fractal set the method consists in hypothetically folding a set of wires according to the geometric shape of each element in the set. The folded wires act as elastic springs. With these springs it is easy to build simple harmonic oscillators. Each virtual harmonic oscillator has clearly three degrees of freedom. Therefore it is possible to find three independent excitations leading to three fundamental periods. The three independent periods provide a complete information set to characterize the topological properties of the sequence. Therefore we may say that the dynamical characterization is more precise than the usual tests as box counting to give an example. The three distinct degrees of freedom generate three fundamental periods which may, or may not, indicate a unique dynamical dimension. For the classical cases of Koch curves where the fractal sequences have well determined Hausdorff dimensions[2] the three distinct dynamical tests lead to a unique fractal dimension.

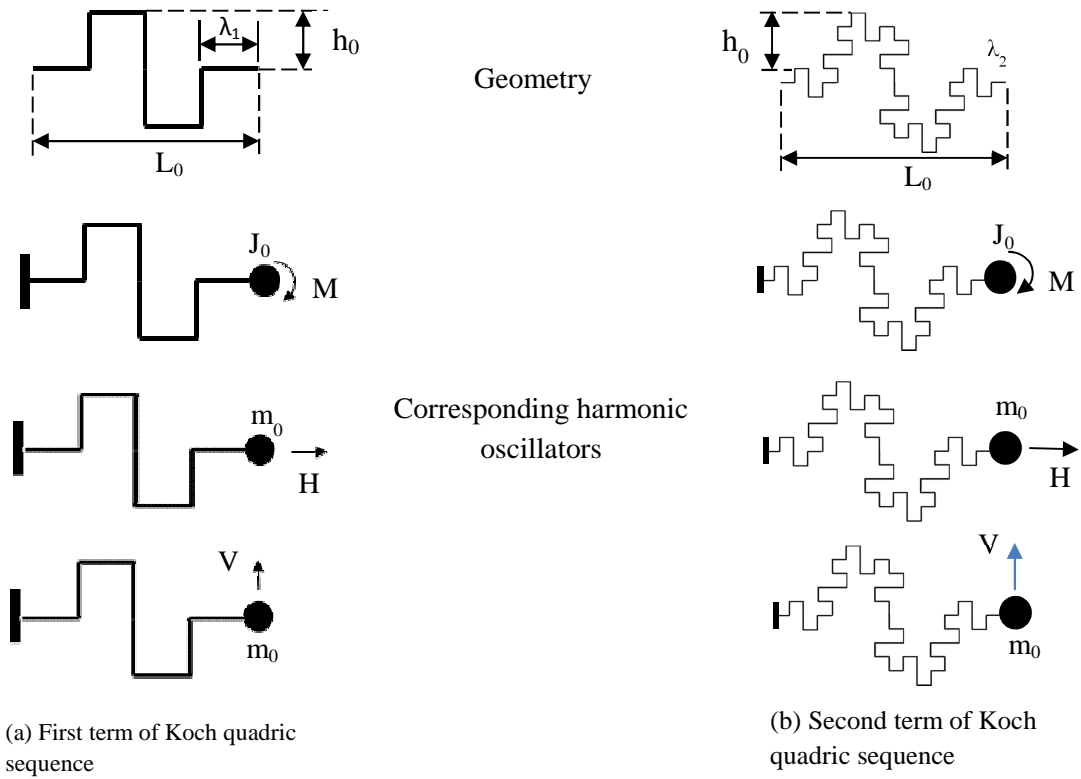

(a) First term of Koch quadric sequence

(b) Second term of Koch quadric sequence

Fig. 1. The Koch quadric curve and the corresponding oscillators with the three distinct excitations characterizing there fundamental frequencies

As shown in Fig.1 we may build a well determined "spring" for each term $k$ of the sequence, folded according to the geometry of the quadric generator, with arbitrary mechanical – Young modulus $E_0$ – and cross section characteristics – moment of inertia $I_0$. Keeping one of the extremities of the spring fixed and attaching to the other a mass $m_0$ we build a harmonic oscillator. The elastic and geometric characteristics $E_0$ and $I_0$ are kept the same for all oscillators. Also the mass $m_0$ attached to the end is kept constant for all elements and the related rotational inertia $J_0$ as well.

Let us impose three different excitation, $M,V,H$ on each element of the sequence. Each of the excitations induces a well determined period that may be written under a normalized form as: $(T_k^M/T_0^M),(T_k^V/T_0^V),(T_k^H/T_0^H)$ for the $k^{th}$ oscillator. The normalizing periods are given by:

$$(T_{M_0})^2 = \frac{J_0 L_0}{E_0 I_0} \qquad (T_{V_0})^2 = \frac{m_0 L_0^3}{E_0 I_0} \qquad (T_{H_0})^2 = \frac{m_0 h_0^2 L_0}{E_0 I_0}$$

The sequences of the normalized periods for the three cases as function of the shortest element $\lambda_k$ represented in a log-log scale approach asymptotically a straight line with a well determined slope that determines the dynamical fractal dimension Fig.2. Let us define the

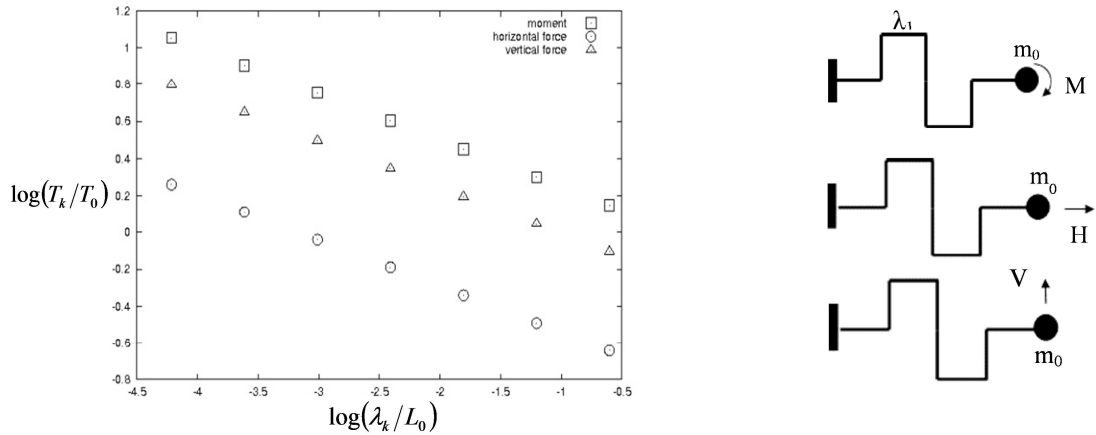

Fig. 2. Koch quadric curve $\log(T_k/T_0)$ x $\log(\lambda_k/L_0)$ for the three excitations M, H, V. $D_M = 1.5$, $D_H = 1.5$  $D_V = 1.5$. Scale of the horizontal projections $L_i/L_0 = (1/1.3)^i$

dimension for this approach as:

$$D_F^\bullet = \lim_{m \to \infty} \left( \frac{\log(T_{m+1}^F / T_m^F)}{\log(\lambda_{m+1}/\lambda_m)} \right)$$

where $F$ represents the dynamical excitations $M,V,H$ respectively.

For the classical Koch curves all three cases provide the same dimension which is related to the Hausdorff geometric dimension $D$ as:

$$D_{M,V,H}^\bullet = \frac{1}{2}(1-D) + \varepsilon$$

where $\varepsilon$ stays for small deviations when the number of elements $m$ is small.
That is, given the first $m$ terms of a fractal sequence, for sufficiently large $m$, the dynamical approach provides an accurate value for the geometric fractal dimension of the sequence. We call the method described above as "direct method".

Now the direct problem applies if we know the elements of the whole sequence or the respective formation rule. The more common case however is the identification problem. That is, given a certain sample, that we will call the master curve, find if it belongs to a fractal sequence and additionally what is the respective fractal dimension. There are several methods easily found in the literature showing successful applications, as the box counting method[2,3]. The most important purpose of those methods is to determine the Hausdorff fractal dimension of the sequence supposedly containing the given sample or the master curve.
Now an alternative technique shows[1,4] that it is also possible to determine the fractal dimension of a "ghost" sequence by cutting successive samples of a given element taken from a "hidden" sequence. The normalized periods $T_n^F/T_0^F$ of each sample $n$ related to its relative projection on the horizontal axis $b_n = L_n/L_0$ represented in a log-log scale provide concrete information about the fractal dimension of some ideal "hidden" sequence. We call it "ghost" sequence because the fractal characteristics obtained with the dynamical approach transfer the actual characteristics of the given sample to all the elements of the hidden sequence. That is, the information contained in the sample is extrapolated to build up the entire hidden sequence. This is not critical for the geometric identification of single fractal structures. However it may be important for the case of multifractals and for lab experiments using a physical sample. The inverse problem formulation cannot detect changes in material and geometrical characteristics of the wire cross section that may change from term to term except if some extra information is given. The fundamental information is related exclusively with the geometry of the curve. Let us define the dynamical fractal dimension for the inverse problem as follows. Consider the set of points $U(\tau_n^F, b_n)$ where $\tau_n^F = T_n^F / T_0^F$ and $b_n = L_n/L_0$. Interpolate a straight line through the points of the set $U$ to obtain:

$$\log(\tau^F) = \overline{D}_F^\bullet \log(b) + \gamma_F$$

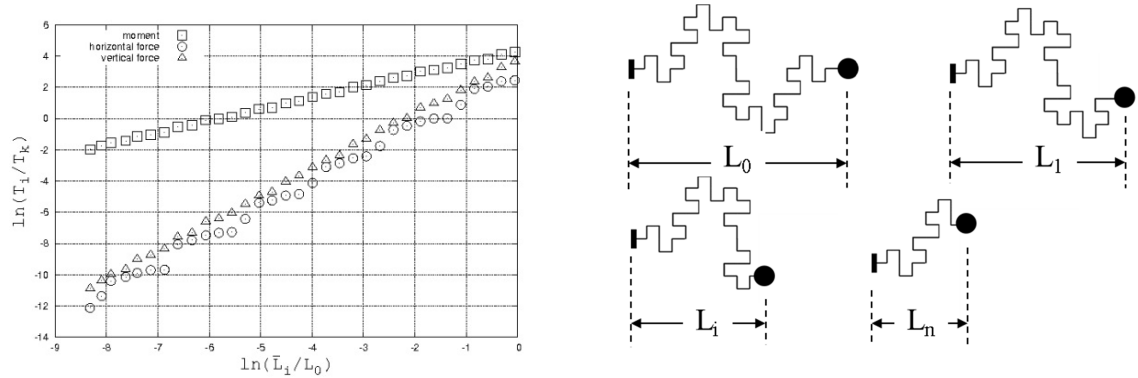

Fig. 3. Koch quadric curve $\log(T_k/T_0) \times \log(\lambda_k/b)$ for the three excitations M, H, V.
$D_M = 1.498505$, $D_H = 1.497955$, $D_V = 1.497938$. Scale of the horizontal projections $L_i/L_0=(1/3)^i$

The slope $\overline{D}_F^\bullet$ is the dynamical fractal dimension corresponding to the excitation $F$ and $\gamma_F$ is the spurious effect introduced by the approximation process.

For some cases like the Koch family curves it is possible to show the correlation between the dynamic dimensions $D_M^\bullet, D_H^\bullet, D_V^\bullet$ and the Hausdorff fractal dimension $D$. Namely:

$$D_M^\bullet = \frac{D}{2} + \varepsilon_M$$

For the case of the initial displacement induced by a moment and

$$D_H^\bullet = D_V^\bullet = \left(1 + \frac{D}{2}\right) + \varepsilon_{H,V}$$

for the case of the initial displacement induced by a vertical or horizontal force. The terms $\varepsilon_M$ and $\varepsilon_{H,V}$ represent the deviations from the true value. These spurious perturbations appear clearly in Fig.3. However the deviations do not induce critical inaccuracies in the determination of the fractal dimension. As shown in Fig. 3, the strongest noise appears in the sequence corresponding to the periods related to the initial condition induced by the horizontal force. The adjusted straight line however gives the same fractal dimension as the other two cases. The maximum absolute error is of the order of 0.013%.

That the dynamical characterization of the "ghost" sequence leads to two different dynamical dimensions is explained by the fact that the energy distribution on the elastic springs due to the bending moment is independent from the position of the segments $\lambda_k$ with respect to a fixed reference frame. For the other two initial perturbations, $H$ and $V$, the energy distribution depends on the position of the elements $\lambda_k$.

Ghost sequences even for classical curves have a two-folded dynamical fractal characterization. This is due to the distinct "energy covers" associated to the imposed excitations. Figure 4 shows the differences among three distinct elastic energy distributions introduced by the bending moments stored in the bars AC-CB. If the acting force is horizontal which is associated with one of the fundamental excitations clearly the stored energy depends on the relative position of the point C with respect to A and B. So given two points A and B and two bars with the same length $\lambda_k$ connecting A and B the stored bending energy is not unique it depends on the position of C. Therefore the horizontal excitation is able to distinguish the position of the bars and consequently the topology of the vertical distribution of the elementary segments of a plane curve. Now for the action of a moment, also associated to one of the fundamental frequencies, the relative positions of the bars is irrelevant. The response to the excitation induced by a moment depends only on the length of the bars. We may say that the "energy cover" corresponding to the moment has similar properties as the Hausdorff cover. It is connected to the topology of the sequence characterized by the union of the elementary segments $\lambda_k$ similar to the Hausdorff cover.

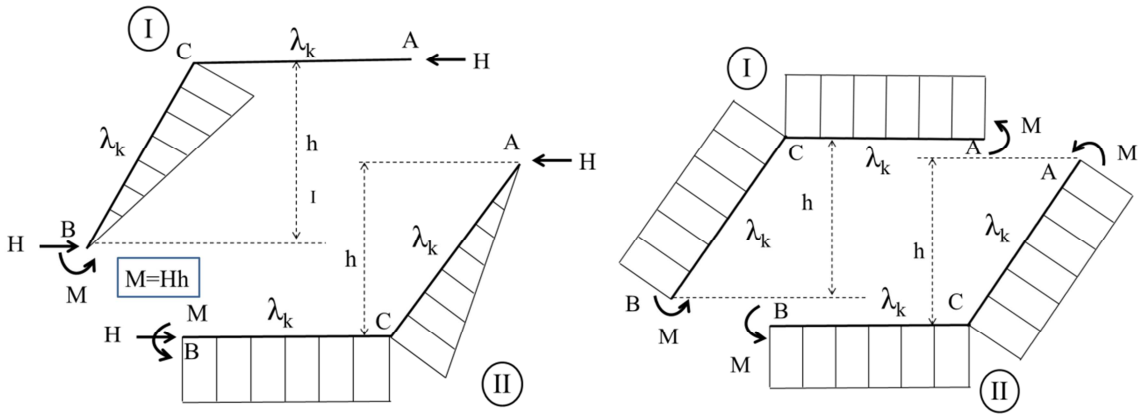

Fig. 4. Bending energy distribution for two configurations I and II due to the action of a horizontal force H and a moment M.

**THE ENERGY COVERING**

The dynamical determination of the fractal dimension is sustained by three different excitations corresponding to three distinct periods. The question that comes out naturally is in what sense are these three sequences different. Does each excitation reveal a different aspect of the fractal characteristic of the curve under examination?

To answer this question let us first recall that the periods are associated to the potential energy stored in the oscillators. Note that the periods are determined using the bending energy only. Therefore if each initial perturbation, *M*, *H* or *V* has a distinct characteristic we may say that each one determines a particular aspect of the "fractality" of the given curve. This is indeed the case. The energy induced by the initial displacement corresponding to a moment M is equally distributed over all the elementary segments $\lambda_k$ for the $k^{th}$ element. Therefore the fractal characterization generated by this initial perturbation depends on the total length of the $k^{th}$ element. Indeed the bending energy stored in this element is:

$$W_k = \frac{1}{2}\frac{1}{EI}\sum_{i=1}^{N_k}\lambda_{i,k} = \frac{1}{2}\frac{1}{EI}L_k$$

where $L_k = N_k \lambda_k$ is the total length of the fractal element. Note that the elementary segments $\lambda_{i,k}$ are all equal $\lambda_{i,k} = \lambda_k$ for the same $k^{th}$ element in the series. The normalized period is easily obtained:

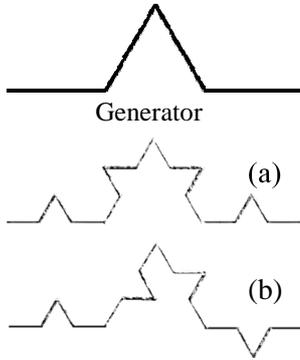

Generator

(a)

(b)

Fig.5. Self-similar grouping (a) and random grouping (b) for the Koch triadic

$$\left(\frac{T_k}{T_0}\right)^2 = \frac{L_k}{L_0} \quad \text{or} \quad \log\left(\frac{T_k}{T_0}\right) = \frac{1}{2}\log\left(\frac{L_k}{L_0}\right)$$

If the curve is fractal with respect to the total length the plot $\log(T_k/T_0) \times \log(\lambda_k/L_0)$ will show this characteristic Fig.2.

Now for the other two cases the bending energy distribution over the elements depends also on the spatial distribution of the elementary segments. The question can now be raised if it is possible, for the dynamical dimension, that a curve presents different dimensions for different excitations, or even if it is possible for a curve to have a fractal characterization for some particular excitation and no fractal characterization for another excitation. In other words is it possible that the dynamical

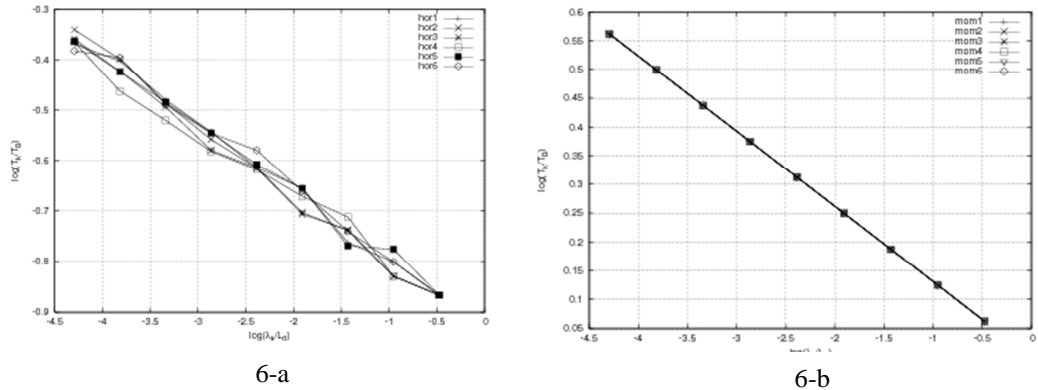

6-a                                                                 6-b

Fig.6. Six realizations of random grouping of Koch triadic generator. The curve $\log(T_k/T_0) \times \log(\lambda_k/b)$ for the case of moment M do not detect the random distribution (1-b). For excitation induced by the horizontal force H (1-a) the curve $\log(T_k/T_0) \times \log(\lambda_k/b)$ clearly shows the lack of self-similarity

dimension depends on the energy covering? This is indeed the case. It was shown[1] that the random orientation of the Koch triadic generator properly downscaled to assemble a triadic like sequence, Fig.6, can be detected by the dynamical response of each term in the sequence. The energy cover generated by the moment doesn't distinguish the self-similar sequence from the random sequence, Fig. 6-b. The energy cover induced by the horizontal force however can clearly show that the assemblage was performed with a random orientation of the Koch triadic generator Fig. 6-a. Clearly the energy cover for the case of the moment is independent of the orientation of the reduced elements in the curve while the orientation of the elements is an important data to evaluate the total energy induced by the horizontal force. For this case due to the strong characterization of the fractal property of the Koch triadic the energy cover derived from the horizontal force reveals fractal characteristics that are distributed in the neighborhood of the self-similar assemblage. The *rms* value calculated from the values given by the different grouping converge to the triadic dimension D=1.26186 .

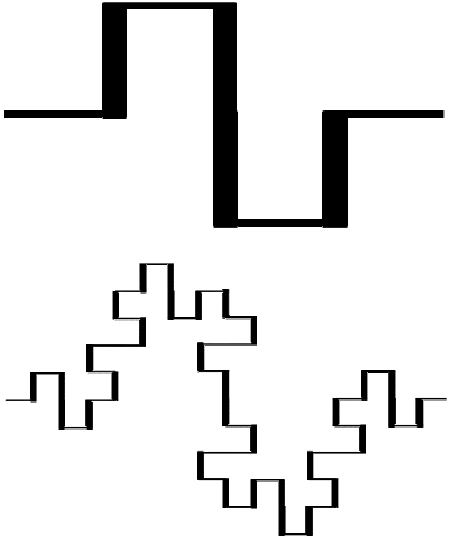

Fig.7. Generator and the first term of the quadric Koch curve. Vertical bars rigid and horizontal flexible

An useful test to verify the characterization of the three different dynamical excitations is to build harmonic oscillators considering the vertical bars rigid and the horizontal bars flexible as shown in Fig. 7. Taking a sufficiently number of terms in the sequence – as a rule the set of the first 12 terms is enough – it is possible to obtain the fractal dimension by plotting the curves:

$$\log\left(\frac{T_k}{T_0}\right)_{M,H,V} = \frac{1}{2}\log\left(f\left(\frac{\lambda_k}{L_0}\right)\right)$$

for the three different excitations *M,H,V*. If the curves come out to be straight lines the slope of these lines will determine the dynamical fractal dimension. For the Koch quadric with rigid vertical bars we obtain:

$D_M = 1.558880$   $D_H = 1.552680$   $D_V = 1.550688$

Given that the geometric fractal dimension is D=1.5 for the flexible configuration the deviations due to introduction of rigid components with the dynamical procedure remain less than 0,4%. The topological identification obtained with the three energy covers clearly reveals a single fractal dimension very close to the case of the flexible configuration. If we invert the flexibility characteristics of the bars, taking the horizontal rigid and the vertical flexible, the result remains unaltered. Therefore we may say that the dynamical method used to evaluate the fractal characteristic of a plane curve uses three different energy covers associated to three distinct topologies. It is important to remark that we are using here the term energy cover in the broad sense. More precisely we may refer to a generalized or extended energy cover meaning that not always there is a classical measure associated to the cover. The energy cover induced by the action of a moment *M* is indeed a cover in the strict sense. It is embedded in a space where a well-defined measure can be introduced. For the cases of the energy distribution induced by a horizontal or vertical force the adoption of the classical measure theory fails. The corresponding embedding space is not assigned with a measure in the classical sense. Therefore when we talk about cover for the energy distribution it has to be understood that we need a new definition, at least for *V* and *H*, where the measure $\mathcal{M}=\mathcal{M}(x_i,y_i,\alpha)$ is function of three independent parameters on a plane. The classical theory has to be extended to include a new parameter.

The above results suggest the introduction of the following propositions:

*i) The energy cover due to the moment excitation is related to the curve length topology.*

*ii) The energy cover due to the horizontal force excitation is related to the topology of the vertical evolution of the curve, that is, the distribution of the elementary segments along the vertical axis.*

*iii) The energy cover due to the vertical force excitation is related to the topology of the horizontal evolution of the curve, that is, the distribution of the elementary segments along the horizontal axis.*

Now, recall that we are dealing with two distinct problems:

1. *Direct problem.* We say that the problem is direct if it is given an ordered sequence of curves following a certain formation law. All curves are defined on a constant interval $(0, L_0)$. The dynamic technique examines the existence of fractal topologies of the given sequence.
2. *Inverse problem.* We say that the problem is inverse if it is given a certain element, of some hidden sequence that may have a fractal structure. The dynamical technique is applied to samples cut off sequentially forming a set of curves $C_1, C_2, ... C_n$ whose projections on the horizontal axis follow a decreasing sequence $L_1 > L_2 > ... > L_n$ starting with the given curve $C_0$ whose projection is $L_0$, the same as the span of the hidden sequence. The dynamical technique is applied to find out a ghost sequence whose topology supposedly fits the topology embedded in the hidden sequence, if any.

Definition:

*A set of plane curves is an entirely fractal set if and only if the three dynamical covers reveals the fractal topology of the set.*

*A given master curve is an entirely fractal curve if and only if the sequence of samples cut off the master curve leads to a ghost sequence that exposes a fractal topology for the three distinct dynamical covers.*

*A set is said to be perfectly fractal if it is entirely fractal and if any of the curves belonging the set is an entirely fractal curve.*

**SELECTED EXAMPLES.**

As introduced in the previous sections the dynamical characterization shows that the Koch quadric is a perfectly fractal curve. For this section we have selected three examples to show the application of the energy cover criterion. One example refers to the Weierstrass-Mandelbrot, Fig.8, curve which is a singular curve generated by an infinity sum[5]:

$$W(t) = \lim_{m \to \infty} \left[ \sum_{n=-m}^{m} \frac{1 - \cos b^n t}{b^{(2-D)n}} \right]$$

Let us consider the Weierstrass-Mandelbrot curve truncated at m=100 with b=1.5 and D=1.5. Figure 8 represents the curve after fifteen interactions[6]. The energy covers for the three excitations *M*, *H* and *V* estimate coherent dynamical fractal dimensions for the direct problem:

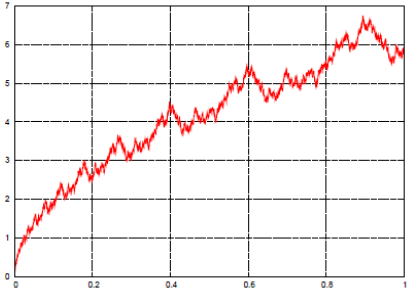

Fig.8. Weierstrass-Mandelbrot curve with D=1.5 for m=100 after 15 interactions.

$D_M=1.463734$; $D_V=1.462550$; $D_V=1,467751$

Note that the fractal dimensions are correlated with the corresponding dynamical fractal dimension as given above:

$$D^{\bullet}_{M_0,V_0,H_0} = \frac{1}{2}(1-D)$$

The errors are less than 2.5%. All the covers indicate a coherent fractal characterization. The sequence representing the Weierstrass-Mandelbrot curve is entirely fractal. Table 1 shows the dynamic dimensions for other values of D obtained with the dynamical approach. The errors remain below ±6%. Taking a large number of terms for the approximating sequences probably the precision could be improved. In any case it is possible to say that all the approximations represent entirely fractal sequences.

|       | D=1.0    | D=1.1    | D=1.3    | D=1.5    | D=1.7    | D=1.9     | D=2.0    |
|-------|----------|----------|----------|----------|----------|-----------|----------|
| $D_M$ | 1.032688 | 1.096462 | 1.288016 | 1.463734 | 1.611834 | 1.7988095 | 2.020544 |
| $D_H$ | 1.044311 | 1.110245 | 1.290945 | 1.462550 | 1.609598 | 1.801783  | 1.998374 |
| $D_V$ | 1.047202 | 1.112634 | 1.293410 | 1.467751 | 1.616178 | 1.808980  | 2.051467 |

Tablel 1. Dynamical dimension for the Weierstrass-Mandelbrot curve evaluated with a converging sequence of curves. Direct problem.

Now if we take the largest term in the Weierstrass-Mandelbrot approximating sequence and try to determine the embedded topology using the sample technique, that is consider the analysis as an inverse problem, the corresponding values obtained for the ghost sequences are shown in the Table 2:

|       | D=1.0    | D=1.1    | D=1.3    | D=1.5    | D=1.7    | D=1.9    | D=2.0    |
|-------|----------|----------|----------|----------|----------|----------|----------|
| $D_M$ | 1.005330 | 1.025217 | 1.031615 | 1.013655 | 1.000221 | 1.000394 | 0.999985 |
| $D_H$ | 1.031976 | 1.151303 | 1.351264 | 1.510396 | 1.675107 | 1.894361 | 2.011184 |
| $D_V$ | 1.010868 | 1.050767 | 1.051465 | 1.011211 | 0.970399 | 0.999697 | 1.008233 |

Table 2. Dynamical dimension for the Weierstrass-Mandelbrot curve using a sequence of samples cutoff from a master curve. Inverse problem.

The dynamical fractal dimension correlates with the fractal dimension as[6]:

$$\overline{D}^{\bullet}_{M_0,V_0,H_0} = 5 - 2D$$

For the inverse problem therefore only the energy cover corresponding to the horizontal force reveals the fractal characteristic of the ghost sequence. Therefore the Weierstrass-Mandelbrot curves do not constitute a perfect fractal set.

Consider now other two cases referring to random generations, namely, random walk and white noise as shown in Fig.9.

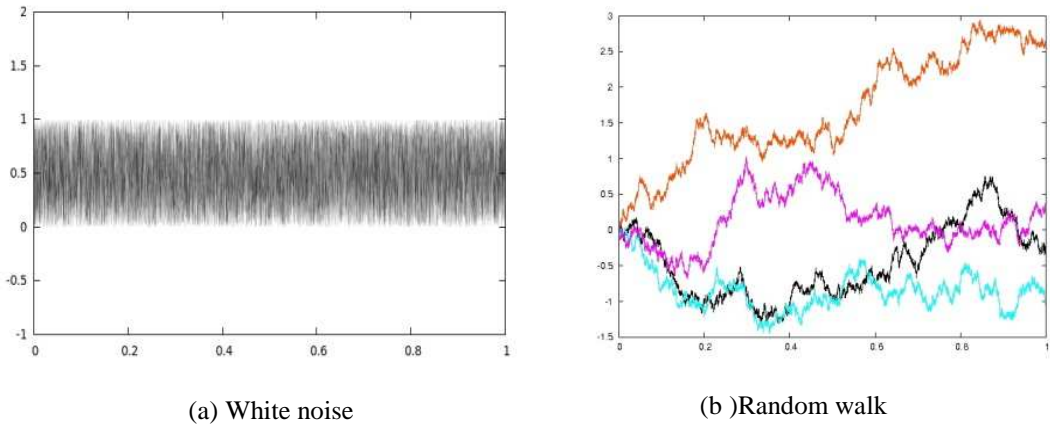

(a) White noise        (b) Random walk

Fig.9. Three particular curves, random walk D=1.5 ; White noise D=2.0 and the Weierstrass-Mandelbrot curve with D=1.5

The random walk is defined in the domain [0,1]. The curve was determined considering $2^{14}$ points and at each step the function deviates of an amount 0.01 upwards or downwards with equal probability. The white noise represents the random distribution in the interval [0,1] long the x-axis of $2^{14}$ points within the interval [0,1] in the y-axis. For these two cases, random walk and white noise we have similar results as for the Weierstrass-Mandelbrot case as shown in Table 3.

| Fractal Dimension | Random walk | | White noise | |
|---|---|---|---|---|
| | Direct problem | Inverse problem | Direct problem | Inverse problem |
| $D_M$ | 1.465275 | 1.000000 | 1.965198 | 1.013887 |
| $D_H$ | 1.498594 | 1.508150 | 1.983766 | 2.033592 |
| $D_V$ | 1.462351 | 0.999995 | 1.981602 | 1.001122 |

Table 3. Dynamical dimension for random walk and white noise, direct and inverse prolems.

As for the Weierstrass-Mandelbrot curve the white noise and the random walk are not perfect fractal sets. For the inverse problem only the energy cover corresponding to the horizontal force characterizes a fractal topology embedded in the master curve.

**CONCLUSION**

The dynamical determination of the fractal topology of plane curves provides more complete information about the character of the fractal structure. It is possible using the dynamic technique to distinguish the topological characteristics of the curve length independently from the distribution of the elementary segments along the horizontal and vertical axis. A curve of the type represented in Fig. 10 is not fractal with respect to the total length although it could present fractal characteristics with respect to the distribution along the vertical axis.

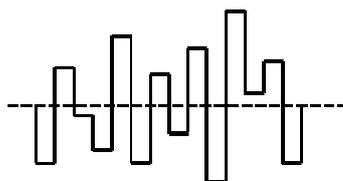

Fig.10. Particular curve with fractal topological distribution along the vertical axis

Note that the mechanical and geometric characteristics of the wires representing the springs of the harmonic oscillators and the mass as well are kept invariants for all terms of the sequences. For more general cases where these characteristics vary for each term of the sequence the fractal characteristic may be considerably modified. Alternating two or more formation rules to build a sequence producing a mixed fractal curve can also be successfully analyzed with the proposed method[7].

The dynamical characterization of fractal curves provides more refined information as compared with other approximation methods to determine the fractal characteristic of plane curves. Since the three distinct energy covering are topologically independent it is possible to detected different fractal formation for particular curves. Affine similarity for instance is a good example where the dynamic dimension may be of great help. Two different similarity scales corresponding to two orthogonal axes may be detected using the dynamic dimension technique. It is also important to mention that the method exposed here may be useful to find the fractal characteristics of fibers and membranes through experiments on material samples. That is the inverse problem may be solved experimentally. Another important application is the determination of the fine fractal structure of proteins. In this case the problem is more complex since we have to deal with harmonic oscillators with at least six degrees of freedom.

The central focus of this paper is to explore possible fractal characteristics of fractal curves using the dynamic response of a set of harmonic oscillators built up according to the geometry of the fractal set. Since the results are promising it is now necessary to explore the deeper analytical meaning underneath the numerical experiments. We believe that new roads concerning the analysis of energy measures and the respective topological consequence should be explored. Also experimental research may be developed to check the efficiency of this method to determine the fractal characteristics of certain physical structures.

## ACKNOWLEDGEMENTS

We are indebted to the National Research Council (CNPq) for the research grants – pos-doc and research scholarships – that made it possible the implementation of the project leading to the results presented here.